\documentclass[12pt]{amsart}
\usepackage{amssymb,amsmath,latexsym,amscd,amsfonts}

\usepackage{pstricks}
\usepackage{multido}
\usepackage[all]{xypic}
\usepackage{ifthen}

\newboolean{showpicture}
\setboolean{showpicture}{true}

\def\gt{\mathfrak{t}}

\def\C{\mathbb{C}}

\def\Z{\mathbb{Z}}

\def\cL{\mathcal{L}}

\def\cP{\mathcal{P}}

\def\cR{\mathcal{R}}

\def\cW{\mathcal{W}}

\def\nek{\text{\hbox{$\simeq$ \kern-.95em \hbox{$/$ \kern.05em}}}}

\def\opp{\operatornamewithlimits{\oplus}}

\def\cplus{\hbox{$\subset${\raise1.05pt\hbox{\kern -0.55em
${\scriptscriptstyle +}$}}\ }}
\def\bcplus{\hbox{$\supset${\raise1.05pt\hbox{\kern -0.55em
${\scriptscriptstyle +}$}}\ }}

\def\ctimes{\hbox{$\times${\raise1.1pt\hbox{\kern -0.27em
${\scriptscriptstyle |}$}}\ }}

\def\udarrow{\hbox{$\nearrow${\kern -0.97em$\searrow$}\ }}

\def\bctimes{\hbox{$\times${\raise1.1pt\hbox{\kern -.74em
${\scriptscriptstyle |}$}}\ }\,\,}

\newcommand{\Agrid}[4]{%
\pscustom[linecolor=gray]{
\code{8 dict begin /lx}
\dim{#1}
\code{def /ly}
\dim{#2}
\code{def /ux}
\dim{#3}
\code{def /uy}
\dim{#4}
\code{def /sc}
\dim{1}
\code{def}
\code{lx ly ux lx sub uy ly sub rectclip}
\code{/horgrid {3 2 roll 3 sqrt div sc div 2 mul floor 2 div sc mul 3 sqrt mul exch 3 sqrt sc mul 2 div exch {3 copy dup 3 1 roll 0 0 [7 1 roll] [1 0 0 1 0 0] 6 array concatmatrix aload pop pop pop moveto lineto pop} for } def }
\code{/lgrid {3 2 roll 3 sqrt div sc div 2 mul floor 2 div sc mul 3 sqrt mul exch 3 sqrt sc mul 2 div exch {3 copy dup 3 1 roll 0 0 [7 1 roll] [60 cos -60 sin  60 sin 60 cos 0 0] 6 array concatmatrix aload pop pop pop moveto lineto pop} for } def }
\code{/rgrid {3 2 roll 3 sqrt div sc div 2 mul floor 2 div sc mul 3 sqrt mul exch 3 sqrt sc mul 2 div exch {3 copy dup 3 1 roll 0 0 [7 1 roll] [-60 cos 60 sin -60 sin -60 cos 0 0] 6 array concatmatrix aload pop pop pop moveto lineto pop} for } def }
\code{lx ly ux uy horgrid}
\code{60 cos ux mul -60 sin ly mul add 60 sin lx mul 60 cos ly mul add 60 cos lx mul -60 sin uy mul add 60 sin ux mul 60 cos uy mul add lgrid}
\code{60 cos lx mul 60 sin ly mul add -60 sin ux mul 60 cos ly mul add 60 cos ux mul 60 sin uy mul add -60 sin lx mul 60 cos uy mul add rgrid} 
\code{end}
}
}

\newcommand{\wtcirc}[2]{\pscircle[fillcolor=gray,fillstyle=solid](!#1 2 div #1 #2 2 mul add 3 0.5 exp mul 6 div){0.1}}

\newcommand{\bigwtcirc}[2]{\pscircle[linestyle=solid](!#1 2 div #1 #2 2 mul add 3 0.5 exp mul 6 div){0.2}} 
\newcommand{\bigwtsq}[2]{
\pscustom{
\msave
\translate(!#1 2 div #1 #2 2 mul add 3 0.5 exp mul 6 div)
\pspolygon(-0.16,-0.16)(-0.16,0.16)(0.16,0.16)(0.16,-0.16)
\translate(!#1 -1 mul 2 div #1 #2 2 mul add 3 0.5 exp mul 6 div -1 mul)
\mrestore
}
}
\newcommand{\wtmvto}[2]{\moveto(!#1 2 div #1 #2 2 mul add 3 0.5 exp mul 6 div)} 
\newcommand{\wtlnto}[2]{\lineto(!#1 2 div #1 #2 2 mul add 3 0.5 exp mul 6 div)}

\newcommand{\Bgrid}[4]{%
\pscustom[linecolor=gray]{
\code{8 dict begin /lx}
\dim{#1}
\code{def /ly}
\dim{#2}
\code{def /ux}
\dim{#3}
\code{def /uy}
\dim{#4}
\code{def /sc}
\dim{1}
\code{def}
\code{lx ly ux lx sub uy ly sub rectclip}
\code{/diaggrid {lx uy sub 2 sqrt div ux ly sub 2 sqrt div lx ly add 2 sqrt div 1 2 sqrt div sc mul ux uy add 2 sqrt div {3 copy dup 3 1 roll 0 0 [7 1 roll] [45 cos -45 sin  45 sin 45 cos 0 0] 6 array concatmatrix aload pop pop pop moveto lineto pop} for } def }  
\code{/odiaggrid {lx ly add 2 sqrt div ux uy add 2 sqrt div lx uy sub 2 sqrt div 1 2 sqrt div sc mul ux ly sub 2 sqrt div {3 copy dup 4 1 roll exch 0 0 [7 1 roll] [45 cos -45 sin  45 sin 45 cos 0 0] 6 array concatmatrix aload pop pop pop moveto lineto pop} for } def }  
\code{ly sc div floor sc mul sc uy sc div ceiling sc mul {dup lx exch moveto ux exch lineto} for}  
\code{lx sc div floor sc mul sc ux sc div ceiling sc mul {dup ly moveto uy lineto} for}            
\code{diaggrid}                                                                        
\code{odiaggrid}                                                                        
\code{end}
}
}

\newcommand{\bwtcirc}[2]{\pscircle[fillcolor=gray,fillstyle=solid](!#1 #2 add #2){0.1}} 
\newcommand{\bigbwtcirc}[2]{\pscircle[linestyle=solid](!#1 #2 add #2){0.2}} 
\newcommand{\bigbwtsq}[2]{
\pscustom{
\msave
\translate(!#1 #2 add #2)
\pspolygon(-0.16,-0.16)(-0.16,0.16)(0.16,0.16)(0.16,-0.16)
\translate(!#1 #2 add -1 mul #2 -1 mul)
\mrestore
}
}
\newcommand{\bwtmvto}[2]{\moveto(!#1 #2 add #2)} 
\newcommand{\bwtlnto}[2]{\lineto(!#1 #2 add #2)} 

\newgray{superlightgray}{0.85}

\newtheorem{theorem}{Theorem}
\newtheorem{example}{Example}
\newtheorem{claim}{Claim}
\newtheorem{question}{Question}
\newtheorem{conjecture}{Conjecture}

\newtheorem{problem}{Problem}

\title[Geometric realization of PRV components]{Geometric realization of PRV components and the Littlewood--Richardson cone}
\author[Ivan Dimitrov]{Ivan Dimitrov$^\dagger$}
\thanks{$^\dagger$ Research partially supported by an NSERC Discovery Grant and by the Max Planck Institute for Mathematics, Bonn}
\author[Mike Roth]{Mike Roth$^*$}
\thanks{$^*$ Research partially supported by an  NSERC Discovery Grant}

\begin{document}


\begin{abstract}
Let $X = G/B$ and let $\cL_1$ and $\cL_2$ be two line bundles on $X$. Consider the
cup product map 
$$
H^{q_1} (X, \cL_1) \otimes H^{q_2} (X, \cL_2) \to H^q(X, \cL),
$$
where $\cL = \cL_1 \otimes \cL_2$ and $q = q_1 + q_2$. We find necessary and sufficient conditions for this map to be
a nonzero map of $G$--modules. We also discuss the converse question, i.e. given irreducible $G$--modules
$U$ and $V$, which irreducible components $W$ of $U \otimes V$ may appear in the right hand side of the equation above.
The answer is surprisingly elegant --- all such $W$ are generalized PRV components of multiplicity one. Along the way we 
encounter numerous connections of our problem with problems coming from Representation Theory, Combinatorics,
and Geometry. Perhaps the most intriguing relations are with questions about the Littlewood--Richardson cone.

This article is expository in nature. We announce results, comment on connections between different fields of Mathematics,
and state a number of open questions. The proofs appear in \cite{dr}.
\end{abstract}

\dedicatory{To Raja on the occasion of his 70$^{\text{th}}$ birthday}

\maketitle


\section*{Introduction}
In 1966 Parthasarathy, Ranga-Rao, and Varadarajan, \cite{PRV}, proved that the tensor product of two irreducible modules $U$ and $V$ 
of a semisimple algebraic group $G$ contains a ``smallest" component $W$ (later named the ``PRV component") whose multiplicity in 
the tensor product is one. The highest weight of $W$ is the dominant weight in the Weyl group orbit of the sum of 
the highest weight of $U$ and the lowest weight of $V$. This remarkable discovery was the first instance of a 
minimal-type representation which later proved to be central in the theory of Harish--Chandra modules.
For details on the history of the PRV component see the excellent article \cite{V}.
In 1988 Kumar, \cite{K1}, generalized the PRV theorem by proving that any irreducible module $G$-module whose 
highest weight is a sum of two extreme weights of $U$ and $V$ is still a component of $U \otimes V$. 
Such components are called ``generalized PRV components". While generalized PRV components 
retain some of the properties of the PRV component, they lack a very important one --- their multiplicities may be greater than one, cf. \cite{K2}. 
This seemed to be the end of the story --- a beautiful discovery with deep applications to Representation Theory, a natural and 
elegant generalization which however lost some features of the original, and it did not seem that there was more to be said.
It turns out, however, that there is another, less obvious generalization of
the PRV component. The construction comes from the cup product on the complete flag variety and yields generalized
PRV components of multiplicity one. 
Moreover, we conjecture that it gives all generalized PRV components of stable multiplicity one. This construction 
gives a different natural generalization of the PRV component --- instead of the ``smallest" component of the tensor
product we obtain ``extreme" components, though not all of them. This last observation leads to natural connections 
with combinatorial problems about the Littlewood--Richardson cone.  

This work arose from the clash between the naively pessimistic intuition derived from Representation Theory of the 
first named author and the optimistic intuition derived from Geometry of the second named author. The truth turned out to
be just a bit off of the latter. The starting point is the Borel--Weil--Bott theorem, \cite{B}, which computes the 
cohomology
of line bundles on complete homogeneous $G$--varieties. Every line bundle has a nonzero cohomology in at most one degree, every such cohomology 
is an irreducible $G$--module, and every 
irreducible $G$--module appears in every degree (not necessarily uniquely) as such a cohomology group. The main application to Representation Theory 
is exactly in constructing all irreducible $G$--modules. In this sense the Borel--Weil theorem, i.e. the statement about cohomology
in degree zero, suffices. As far as we know, Bott's theorem ---the statement about higher cohomology --- has not been
used for constructing representations of reductive algebraic groups in characteristic zero. 
In this paper we apply Bott's theorem to construct irreducible components of the tensor
product of two irreducible representations. Namely, we consider the diagonal embedding of the homogeneous variety
$X=G/B$ , where $B$ is a Borel subgroup of $G$ into $X \times X$. It gives rise to a map $\pi$ from the cohomology
of line bundles on $X \times X$ to the cohomology of the restrictions of these line bundles on $X$. Since the diagonal
embedding of $X$ into $X \times X$ is $G$--invariant, the map $\pi$ is a $G$--module map from the tensor product
$U^* \otimes V^*$ of two irreducible
$G$--modules $U^*$ and $V^*$ to another irreducible $G$--module $W^*$. If the map is nonzero then by dualizing we obtain a geometric construction of 
the simple component $W$ of $U \otimes V$. Two natural problems arise from this situation --- find necessary and sufficient conditions
for $\pi$ to be a nonzero map between $G$--modules and describe all components of $U \otimes V$ which can be constructed in this way.

We solve the first problem by finding an explicit (though somewhat mysterious) necessary and sufficient combinatorial condition on the
line bundles under consideration. The condition, equation (\ref{equation6}), is expressed in terms of a triple of Weyl group elements
naturally associated to the line bundle on $X \times X$. 
The second question seems to be more difficult and we only have a partial solution and a conjecture about the full answer. 
We prove that the geometric construction yields only generalized PRV components of multiplicity one and we conjecture
that every generalized PRV component of stable multiplicity one can be obtained in this way. Furthermore, we show that
the components that result from the geometric construction are always extreme in the Littlewood--Richardson cone. This
relates our results to the recent renewed interest in the structure of the Littlewood--Richardson cone which began with 
the proof of Horn's conjecture in the work of A. Klyachko, \cite{Kly}, and A. Knutson and T. Tao, \cite{KT1}. 

The present paper is expository in nature. We state the two problems discussed above and the main theorems that
relate to them as well several other statements of interest.  We also pose a number of open questions. 
The proofs appear in \cite{dr}. 
Here is briefly the content of each of the five sections.

\begin{itemize}
\item[1.] Statement of the two main problems.
\item[2.] Solution of the first problem.
\item[3.] Discussion of the second problem.
\item[4.] Relation of the second problem for $GL_{n+1}$ to Schubert calculus.
\item[5.] Relation of our results to the structure of the Littlewood--Richardson cone and to other ``cone" problems.
\end{itemize}

\noindent
{\bf Acknowledgments.} We thank P. Belkale, W. Fulton, B. Kostant, S. Kumar, K. Purbhoo, and N. Reading for many fruitful discussions.
I. D. acknowledges the support and excellent working conditions at the Max Planck Institute, Bonn. M. R. 
acknowledges the hospitality of the University of Roma III.

\medskip
\noindent
{\bf Notation.} The ground field is $\C$.
Let $G$ be a connected reductive algebraic group, let $T \subset B \subset G$ be a maximal torus and a Borel subgroup of $G$,
and let $X = G/B$. We have the following standard objects related to the triple $T \subset B \subset G$.
\begin{itemize}
\item $\cP$ --- the characters of $T$ and $\cP^+$ --- the dominant characters of $T$;
\item $V(\lambda')$ --- the irreducible $G$--module with $B$--highest weight $\lambda' \in \cP^+$;
\item $\Delta$ --- the roots of $G$, $\Delta^+$ --- the roots of $B$, and $\Delta^- := - \Delta^+$;
\item $\cW$ --- the Weyl group of $G$. For $w \in \cW$ we denote the length of $w$ by $l(w)$;
\item $w_0 \in \cW$ --- the longest element of $\cW$;
\item $\rho := 1/2 \sum_{\alpha \in \Delta^+} \alpha$;
\item $(-, -): \gt^* \times \gt^* \to \C$ --- a non--degenerate $\cW$--invariant symmetric bilinear form on $\gt^*$, where $\gt := {\rm{Lie }}\, \, T$;
\item $c_{\lambda', \mu'}^{\nu'}$ --- the multiplicity $[V(\lambda') \otimes V(\mu') : V(\nu')]$.
\end{itemize}
We consider two actions of $\cW$ on $\cP$ --- the usual action which we call homogeneous and denote by $w \lambda$ or $w(\lambda)$,
and the affine action given by $w \cdot \lambda := w (\lambda + \rho) - \rho$. 
A character $\lambda \in \cP$ is regular if there is a dominant character in the affine orbit of $\lambda$, or equivalently if 
$(\lambda + \rho, \alpha) \neq 0$ for every $\alpha \in \Delta$. If $\lambda$ is regular, then there exists a unique element $w_\lambda \in \cW$
for which $w_\lambda \cdot \lambda$ is dominant. We define the length of $\lambda$, $l(\lambda)$, to be the length of $w_\lambda$. If
$\lambda$ is not regular, we call it singular. The length of a singular element is not defined. Whenever we use the notation $l(\lambda)$ we
assume implicitly that $\lambda$ is regular.

\section{The Borel--Weil--Bott theorem and diagonal embeddings.}
For $\lambda \in \cP$, let $\cL_\lambda$ be the line bundle on $X$ corresponding to the $B$--module $\C_{-\lambda}$ on which
$T$ acts via the character $(-\lambda)$ and the unipotent radical of $B$ acts trivially. The Borel--Weil--Bott theorem, see \cite{B}, states that
$$
H^q(X, \cL_\lambda) = \left\{ \begin{array}{cl}
V(w_\lambda \cdot \lambda)^* & {\text { if }} l(\lambda) = q\\
0 & {\text { otherwise}}. 
\end{array} \right.
$$

Consider the diagonal embedding $X \hookrightarrow X \times X$. If $\lambda, \mu \in \cP$, the line bundle 
$\cL_\lambda \boxtimes \cL_\mu$ on $X \times X$ restricts to the line bundle $\cL_{\lambda+ \mu}$ on $X$
and this restriction induces a natural map
\begin{equation} \label{equation2}
\pi: H^q(X \times X, \cL_\lambda \boxtimes \cL_\mu) \to H^q(X, \cL_{\lambda+ \mu}).
\end{equation}
Both sides of (\ref{equation2}) are $G$--modules and the map $\pi$ is a $G$--module homomorphism. If one or both of these
modules are zero, then $\pi$ is trivial. Assume both sides of (\ref{equation2}) are nontrivial $G$ modules, i.e. assume that
\begin{equation} \label{equation3}
l(\lambda + \mu) = l(\lambda) + l(\mu),
\end{equation}
and that $q$ is the common value. The Borel--Weil--Bott theorem allows us to compute explicitly the $G$--modules in (\ref{equation2}).
By Kunneth's theorem we have
\begin{align*}
& H^q(X \times X, \cL_\lambda \boxtimes \cL_\mu) = \opp_{i+j = q} H^i(X, \cL_\lambda) \otimes H^j(X, \cL_\mu) \\ = &
H^{l(\lambda)}(X, \cL_\lambda) \otimes H^{l(\mu)}(X, \cL_\mu) 
 =  V(w_\lambda \cdot \lambda)^* \otimes V(w_\mu \cdot \mu)^*.
\end{align*}
Hence, assuming (\ref{equation3}), the dual of $\pi$ is a $G$--module homomorphism
\begin{equation} \label{equation4}
\pi^*: V(\nu') \to V(\lambda') \otimes V(\mu'),
\end{equation}
where $\lambda' := w_\lambda \cdot \lambda$, $\mu':= w_\mu \cdot \mu$, and $\nu' := w_{\lambda + \mu} \cdot (\lambda + \mu)$.
Since $V(\nu')$ is an irreducible $G$--module, $\pi^*$ is zero or injective and, respectively, $\pi$ is zero or surjective.
Thus we have arrived at the two main problems of the present paper.

\begin{problem} \label{problem1}
When is $\pi$ a surjective map between nontrivial $G$--modules?
\end{problem}

\begin{problem} \label{problem2}
Given $\lambda', \mu' \in \cP^+$, find all simple components $V(\nu')$ of $V(\lambda') \otimes V(\mu')$ which arise from (\ref{equation4}).
More precisely, find all $\nu' \in \cP^+$ for which there exist $w_1, w_2, w_3 \in W$ with the property that 
\begin{equation} \label{equation22}
w_3^{-1} \cdot \nu' = w_1^{-1} \cdot \lambda' + w_2^{-1} \cdot \mu'
\end{equation}
 and such that the map $\pi^*$ corresponding to
$\lambda = w_1^{-1} \cdot \lambda'$ and $\mu = w_2^{-1} \cdot \mu'$ is an injective map between nontrivial 
$G$--modules.\footnote{\, Strictly speaking, we need to do more if the multiplicity of $V(\nu')$ in $V(\lambda') \otimes V(\mu')$ is greater than one.
Since this never happens, describing all characters $\nu'$ suffices.} 
\end{problem}
We call such components $V(\nu')$ of $V(\lambda') \otimes V(\mu')$ cohomological.

\section{Inversion sets and the answer to Problem 1.}
We start our discussion with the first problem. Clearly, (\ref{equation3}) is a necessary condition. Unfortunately, it is not sufficient
as the following example shows.

\begin{example} \label{example1} Let $G = GL_6$ and let $\lambda = \mu = (a, b, 0, a+2, b+2, 2)$ for some integers $a > b > 0$.
 One checks immediately that $l(\lambda) = 3$, while $l(2 \lambda) = 6$.
Furthermore, $\lambda' = \mu' = (a, a, b+1, b+1, 2, 2)$ and $\nu' = (2a+1, 2a+1, 2b+2, 2b+2, 3,3)$, which shows that $V(\nu')$ is 
not a component of $V(\lambda') \otimes V(\mu')$, and hence $\pi$ and $\pi^*$ are trivial.
\end{example}

To state the answer to Problem 1 we need to introduce some more notation. The inversion set $\Phi_w$ of an element $w \in \cW$
is defined as the set of positive roots which $w$ sends to negative roots, i.e. $\Phi_w := \Delta^+ \cap w^{-1} \Delta^-$. Set 
$\Phi_w^c := \Delta^+ \backslash \Phi_w$. Inversion sets were introduced by Kostant in \cite{Ko}. 
Assume now that for $\lambda, \mu \in \cP$ condition (\ref{equation3}) holds. Note first that 
\begin{align*}
l(\lambda) & = \# \{ \alpha \in \Delta^+ \, | \, (\lambda + \rho, \alpha) < 0 \} 
 = \# \{ \alpha \in \Delta^+ \ | \, (w_1(\lambda + \rho), w_1 \alpha) <0\} \\ &
= \# \{ \alpha \in \Delta^+ \, | \, (\lambda' + \rho, w_1 \alpha ) < 0 \} = \# \Phi_{w_\lambda},
\end{align*}
where $\# S$ stands for the cardinality of a set $S$. Similarly, $l(\mu) = \# \Phi_{w_\mu}$ and $l(\lambda + \mu) = \# \Phi_{w_{\lambda + \mu}}$.
Condition (\ref{equation3}) is therefore equivalent to 
\begin{equation} \label{equation5}
\# \Phi_{w_{\lambda + \mu}} = \# \Phi_{w_\lambda} + \# \Phi_{w_\mu}.
\end{equation}
Assume additionally that $\lambda + \rho$ and $\mu + \rho$ are sufficiently far from the walls of the Weyl chambers. (Explicitly, it is enough to assume that
for every $\alpha \in \Delta$, $|(\lambda + \rho, \alpha)|$ and $|(\mu + \rho, \alpha)|$ are greater than $1/2$ the maximal height of a root of $G$.)
Then $\Phi_{w_{\lambda + \mu}} \subset \Phi_{w_\lambda} \cup \Phi_{w_\mu}$. Indeed, if $\alpha \in \Phi_{w_\lambda}^c \cap \Phi_{w_\mu}^c$, then
$$
(\lambda + \mu + \rho, \alpha) = (\lambda + \mu + 2 \rho, \alpha) - (\rho, \alpha) = (\lambda + \rho, \alpha) + 
(\mu + \rho, \alpha) - (\rho, \alpha) > 0, 
$$  
i.e. $\alpha \in \Phi_{w_{\lambda+\mu}}^c$. Now $\Phi_{w_{\lambda+\mu}} \subset \Phi_{w_\lambda} \cup \Phi_{w_\mu}$ together with (\ref{equation5}) implies 
\begin{equation} \label{equation6}
\Phi_{w_{\lambda+\mu}} = \Phi_{w_\lambda} \sqcup \Phi_{w_\mu}.
\end{equation}
To summarize the discussion above, we have shown that for $\lambda + \rho$ and $\mu + \rho$ far enough from the walls of the Weyl chambers,
conditions (\ref{equation3}) and (\ref{equation6}) are equivalent. 
Example \ref{example1}, where $\Phi_{w_\lambda}=\Phi_{w_\mu}$, shows that (\ref{equation3}) does not imply (\ref{equation6}) 
in general.
Condition (\ref{equation6}) above is somewhat mysterious,
however it does appear in other instances. Combinatorially it can be expressed in terms of the weak Bruhat order on $\cW$ --- it means that
the greatest lower bound of $w_\lambda$ and $w_\mu$ is the identity and the least upper bound of the two is $w_{\lambda+ \mu}$.\footnote{\, We thank N. Reading for telling
us about this interpretation.} Another interesting fact is that (\ref{equation6}) is equivalent to the property that $s(\Phi_{w_\lambda}), s(\Phi_{w_\mu})$,
and $s(\Phi_{w_{\lambda + \mu}}^c)$ are mutually orthogonal, where $s(\Phi) = \sum_{\alpha \in \Phi} \alpha$ for any subset $\Phi \subset \Delta^+$.\footnote{\, We thank 
B. Kostant for telling us this fact.} The appearance of (\ref{equation6}) most relevant to our work is in \cite{BK1} where
Belkale and Kumar define a new product in the ring $H^*(X, \Z)$ by keeping the structural constants 
$d_{w_1, w_2}^{w_3}$ (see (\ref{equation11}) below) corresponding to triples satisfying $\Phi_{w_3} = \Phi_{w_1} \sqcup \Phi_{w_2}$ the same, and setting all other 
structure constants equal to zero. This alternate product on $H^*(X, \Z)$ is then used to parameterize a minimal set of inequalities determining the
cone of solutions of an eigenvalue problem associated to $G$.

The solution of Problem 1 is given by the following theorem.

\begin{theorem} \label{theorem1} 
The map $\pi$ is a surjection of nontrivial $G$--modules if and only if 
$\Phi_{w_{\lambda + \mu}} = \Phi_{w_\lambda} \sqcup \Phi_{w_\mu}$
and 
$q = l(\lambda + \mu)$.
\end{theorem}

Here are a few words about the proof. The cohomology ring $H^*(X, \Z)$ plays a crucial role. For each $w \in \cW$, let
$X_w :=\overline{BwB/B} \subset X$ denote the Schubert variety associated to $w$.
The classes $\{[X_w]\}_{w \in \cW}$ of all
Schubert varieties form a basis of $H^*(X, \Z)$.  The Poincar\'e dual basis $\{[\Omega_w]\}_{w \in \cW}$ is given by $\Omega_w := X_{w_0 w}$, 
where $w_0$ denotes the longest element of $\cW$, cf. \cite{De1}.   If $l(w_3) = l(w_1) + l(w_2)$ set
\begin{equation} \label{equation11}
d_{w_1, w_2}^{w_3} := ([\Omega_{w_1}] \cap [\Omega_{w_2}]) \cdot [X_{w_3}].
\end{equation}
If we assume that $\pi$ is a surjection of nontrivial $G$--modules, then we first show that $d_{w_\lambda, w_\mu}^{w_{\lambda+\mu}} \neq 0$.
Furthermore, $d_{w_\lambda, w_\mu}^{w_{\lambda+\mu}} \neq 0$ together with $c_{\lambda', \mu'}^{\nu'} \neq 0$ 
implies that 
$\Phi_{w_{\lambda + \mu}} = \Phi_{w_\lambda} \sqcup \Phi_{w_\mu}$ which establishes one direction of
the theorem. The other direction is more difficult to prove. The first step is to see that $d_{w_\lambda, w_\mu}^{w_{\lambda + \mu}} \neq 0$ and then
a rather delicate geometric argument completes the proof. This last step is somewhat simpler if $d_{w_\lambda, w_\mu}^{w_{\lambda + \mu}} = 1$.
In fact it seems that this may always be the case.

\begin{claim} \label{claim0}
If $G$ is a simple classical group, then
$\Phi_{w_3} = \Phi_{w_1} \sqcup \Phi_{w_2}$ implies that $d_{w_1, w_2}^{w_3} = 1$.
\end{claim}
P. Belkale and S. Kumar showed us a proof of Claim \ref{claim0} in the case of $G = SL_{n+1}$, \cite{BK2}. Their proof goes through for simple
groups of type $B$ and $C$ as well. The case of simple groups of type $D$ is more difficult and involves both combinatorial and
geometric machinery. It is also not difficult to check the statement above for the group $G_2$. 
The remaining exceptional groups (at least $E_8)$ seem to be beyond a computer verification. 

We complete this section by stating two open questions.

\begin{question} \label{question0}
Is it true that $\Phi_{w_3} = \Phi_{w_1} \sqcup \Phi_{w_2}$ implies that 
$d_{w_1, w_2}^{w_3} = 1$ for any $G$?
\end{question}

\begin{question} \label{question22}
Is it true that  if $l(\nu)=l(\lambda) + l(\mu)$ then $\pi$ is surjective map of nontrivial $G$--modules if and only if $c_{\lambda',\mu'}^{\nu'} \neq 0$?
\end{question}

\section{Generalized PRV components and cohomological components.}
In this section we discuss Problem \ref{problem2}. If $\lambda', \mu' \in \cP^+$ and $\nu'$ is the dominant character in the
$\cW$-orbit of $\lambda'+w_0\mu'$,
where $w_0$ is the longest element of $\cW$, then $V(\nu')$ is a component of $V(\lambda') \otimes V(\mu')$ of multiplicity one, see \cite{PRV} . 
The component $V(\nu')$ is called the PRV component of $V(\lambda') \otimes V(\mu')$. More generally, for any $w \in \cW$, 
if $\nu'$  is the dominant character in the $\cW$-orbit of $\lambda'+w\mu'$ then
$V(\nu')$  is still
a component of $V(\lambda') \otimes V(\mu')$, see \cite{K1}. These are called generalized PRV components of $V(\lambda') \otimes V(\mu')$. 
Unlike the PRV component, generalized PRV components may have multiplicities greater than one. We first show that every cohomological
component of $V(\lambda') \otimes V(\mu')$ is a generalized PRV component. Equation (\ref{equation22}) can be rewritten as 
$\nu' = w_3^{-1} \cdot (w_1^{-1} \cdot \lambda' + w_2^{-1} \cdot \mu')$ which resembles the expression relating the highest weights of the 
generalized PRV components except for the fact that we have the affine action of the Weyl group instead of the homogeneous one.
However, if $V(\nu')$ is a cohomological component then Theorem \ref{theorem1} implies that $\Phi_{w_3} = \Phi_{w_1} \sqcup \Phi_{w_2}$.
Using the equality $w^{-1} \cdot 0  = w^{-1} \rho - \rho = - s(\Phi_w)$, where $s(\Phi_w) = \sum_{\alpha\in \Phi_w} \alpha$ as above, we conclude
that 
$\Phi_{w_3} = \Phi_{w_1} \sqcup \Phi_{w_2}$
implies that $w_3^{-1} \rho + \rho = w_1^{-1} \rho + w_2^{-1} \rho$. The last equation ensures that 
$w_3^{-1} \cdot \nu' = w_1^{-1} \cdot \lambda' + w_2^{-1} \cdot \mu'$ is equivalent to
$\nu' = w_3(w_1^{-1} \lambda' + w_2^{-1}  \mu')$, which shows that every cohomological component of $V(\lambda') \otimes V(\mu')$ is
a generalized PRV component. 
Notice that the argument above relies on the fact that 
$$\Phi_{w_3} = \Phi_{w_1} \sqcup \Phi_{w_2}\,\, \mbox{implies that}\,\,
w_3^{-1} \cdot 0 = w_1^{-1} \cdot 0 + w_2^{-1} \cdot 0.
$$
The converse is also true when $d_{w_1, w_2}^{w_3}\neq 0$ as the following statement shows. 

\begin{claim} \label{claim10}
Let $w_1, w_2, w_3 \in \cW$ be such that $l(w_3) = l(w_1) + l(w_2)$ and $d_{w_1, w_2}^{w_3}\neq 0$.
Then $w_3^{-1} \cdot 0 = w_1^{-1} \cdot 0 + w_2^{-1} \cdot 0$ implies $\Phi_{w_3} = \Phi_{w_1} \sqcup \Phi_{w_2}$.
\end{claim}

A partial solution to Problem \ref{problem2} is given by the next theorem.

\begin{theorem} \label{theorem2}
If   $V(\nu')$ is a cohomological component of $V(\lambda') \otimes V(\mu')$ then $V(\nu')$ is a generalized
PRV component of $V(\lambda') \otimes V(\mu')$ of multiplicity one.
\end{theorem}

Theorem \ref{theorem1} implies that if $V(\nu')$ is a cohomological component of $V(\lambda') \otimes V(\mu')$ then
$V(k \nu')$ is a cohomological component of $V(k \lambda') \otimes V(k \mu')$ for every positive integer $k$. This,
combined with Theorem \ref{theorem2}, implies the following.

\begin{claim} \label{claim1}
If $V(\nu')$ is a cohomological component of $V(\lambda') \otimes V(\mu')$ then $c_{k \lambda', k \mu'}^{k \nu'}= 1$ for
every positive integer $k$.
\end{claim}

We believe that the converse of Claim \ref{claim1} is also correct.

\begin{conjecture} $V(\nu')$ is a cohomological component of $V(\lambda') \otimes V(\mu')$ if and only if $V(\nu')$ is a generalized
PRV component of $V(\lambda') \otimes V(\mu')$ of stable multiplicity one, i.e.  $c_{k \lambda', k \mu'}^{k \nu'}= 1$ for every positive integer $k$.
\end{conjecture}

This conjecture is supported by the following particular cases.

\begin{claim} \label{claim2} In each of the following cases, $V(\nu')$ is a generalized PRV component of stable multiplicity one in $V(\lambda') \otimes V(\mu')$.
Moreover, $V(\nu')$ is a cohomological component of  $V(\lambda') \otimes V(\mu')$.
\begin{itemize}
\item[{\rm(i)}] When $V(\nu')$ is the PRV component of $V(\lambda') \otimes V(\mu')$.
\item[{\rm(ii)}] When $\lambda' \gg \mu'$, in the sense that  $\lambda' + w \mu' \in \cP^+$ for every $w \in \cW$, and 
$\nu' = \lambda' + w \mu'$ for some $w \in \cW$.
\end{itemize}
\end{claim}
The fact that $V(\nu')$ above is of stable multiplicity one in $V(\lambda') \otimes V(\mu')$ follows from \cite{PRV} in the first case, and is
an elementary exercise in the second one.
The construction of the corresponding triple $(w_1, w_2, w_3)$ is straightforward in both cases. 
For the PRV component
the triple is given by $(w_0 \sigma^{-1}, \sigma^{-1}, w_0)$, where
$\sigma\in\cW$ is an element so that $\nu'=\sigma(\lambda'+w_0\mu')\in\cP^{+}$. 
In the second case we can simply
take the triple to be $(w^{-1}, e, w^{-1})$. 

For $G = GL_{n+1}$ a conjecture of Fulton, proved by Knutson, Tao, and Woodward \cite{KTW}, 
states that $c_{\lambda', \mu'}^{\nu'} = 1$ is equivalent to  
$c_{k \lambda', k \mu'}^{k \nu'}= 1$ for every positive integer $k$. 
When $G$ is of type $A$ the conjecture would therefore imply that a component $V(\nu')$ of $V(\lambda')\otimes V(\mu')$ is cohomological
if and only if $V(\nu')$ is a generalized PRV component of multiplicity one.
The fact that multiplicity one implies stable multiplicity one is 
not true in general. The next example 
illustrates this and our conjecture for $G = SO_5$.

\begin{example} \label{example2} 
Let $G = SO_5$ and let $\lambda' = \mu' = \rho = \omega_1 + \omega_2$, where $\omega_1$
and $\omega_2$ are the fundamental weights. $V(\rho) \otimes V(\rho)$ contains the following components 
$$
(0,0), (1,0), (2,0), (3,0), (0,2), (1,2), (2,2), (0,4),
$$ 
where $(a,b)$ denotes $a \omega_1 + b \omega_2$ (see the middle picture in Figure 1 for this decomposition). 
The generalized PRV components are 
$$
(0,0), (1,0), (3,0), (0,2), (1,2), (2,2), (0,4)
$$ 
of which $(1,2)$ and $(0,2)$ have multiplicity $2$ and the rest have multiplicity $1$. The cohomological components are
$$
(0,0), (3,0), (2,2), (0,4).
$$ 
The component $(1,0)$ is a generalized PRV component of multiplicity $1$ which is not cohomological. This does not contradict the Conjecture
since for $k=2$ we have $c_{(2,2),(2,2)}^{(2,0)} = 2$.
\end{example}
It is interesting to know for which $k$ we need to check the multiplicities  $c_{k \lambda', k \mu'}^{k \nu'}$. 
Kapovich and Millson, 
\cite{KM}, proved that there exists $k = k(G)$ such that $c_{k \lambda', k \mu'}^{k \nu'} \neq 0$ if and only if 
$c_{N \lambda', N \mu'}^{N \nu'} \neq 0$ for $N \geq k$. However we do not know whether a similar result holds for detecting
stable multiplicity one.

\section{Reduction patterns and the proof of Theorem \ref{theorem2} for $G= GL_{n+1}$.}
The proof of Theorem \ref{theorem2} is type--independent.\footnote{\, As Raja has taught I. D., one does not really understand a theorem about
semisimple groups until there is a type--independent proof.}  However, if $G = GL_{n+1}$ there is a different proof which exploits
the fact that the Littlewood--Richardson coefficients appear as structure constants for the multiplication in the 
cohomology ring of
Grassmannians. Since this proof establishes yet another connection between cohomological components and 
classical geometric objects, we will outline it here.

Let $G=GL_{n+1}$ and $\tilde{G}=GL_{n}$.
Assuming that $V(\nu')$ is a cohomological
component of $V(\lambda') \otimes V(\mu')$ for the group $G$ we will construct a triple $(\tilde{\lambda}', \tilde{\mu}', \tilde{\nu}')$
such that $V(\tilde{\nu}')$ is a cohomological component of $V(\tilde{\lambda}') \otimes V(\tilde{\mu}')$ for the group $\tilde{G}$
and with the property that $c_{\lambda', \mu'}^{\nu'} = c_{\tilde{\lambda}', \tilde{\mu}'}^{\tilde{\nu}'}$. We call the assignment 
$$
\cR_n\colon
(\lambda', \mu', \nu') \to (\tilde{\lambda}', \tilde{\mu}', \tilde{\nu}')
$$
a reduction of $(\lambda', \mu', \nu')$. Since $V(\tilde{\nu}')$ is again a cohomological component of $V(\tilde{\lambda}') \otimes V(\tilde{\mu}')$
there is a reduction $\cR_{n-1}$ of  $(\tilde{\lambda}', \tilde{\mu}', \tilde{\nu}')$ which again preserves the corresponding Littlewood--Richardson coefficients. 
Continuing inductively we obtain a reduction pattern $\cR$, i.e. a composition $\cR = \cR_1 \circ \cR_2 \circ \ldots \circ \cR_n$ 
of consecutive reductions that can be applied to $(\lambda', \mu',\nu')$ which at each step preserves 
the property of being a cohomological component 
and 
preserves the corresponding Littlewood--Richardson coefficient.
Since the tensor product of irreducible $GL(1)$--modules is irreducible we conclude that $c_{\lambda', \mu'}^{\nu'} = 1$.

To find characters $\tilde{\lambda}'$, $\tilde{\mu}'$, and $\tilde{\nu}'$ such that 
$c_{\lambda', \mu'}^{\nu'} = c_{\tilde{\lambda}', \tilde{\mu}'}^{\tilde{\nu}'}$
we will make use of the fact that the Littlewood-Richardson 
coefficients coincide with the intersection numbers of Schubert cycles on Grassmannians.
More precisely, in the notation of \cite{GH}, $c_{\lambda', \mu'}^{\nu'} = \#(\sigma_{\lambda'} \cdot \sigma_{\mu'} \cdot \sigma_{(\nu')^*})$,
where $\sigma_{(\nu')^*}$ is the cycle which is Poincar\'e dual to $\sigma_{\nu'}$.  
The reduction $\cR_{n}$ will be defined as the counterpart of Reduction Formula I on p. 202 of \cite{GH}, i.e. through the diagram
$$
\xymatrix { c_{\lambda', \mu'}^{\nu'} \ar@{=}[r] \ar[d]_{\cR_{n}} &  \#(\sigma_{\lambda'} \cdot \sigma_{\mu'} \cdot \sigma_{(\nu')^*}) \ar[d]^{\mbox{\tiny Reduction Formula I}}\\
c_{\tilde{\lambda}', \tilde{\mu}'}^{\tilde{\nu}'} & \#(\sigma_{\tilde{\lambda}'} \cdot \sigma_{\tilde{\mu}'} \cdot \sigma_{(\tilde{\nu}')^*}). \ar@{=}[l]
}
$$
Here are the explicit formulas for $\cR_{n}$.
Let $\lambda' = (\lambda_0', \lambda_1', \ldots, \lambda_n')$, 
$\mu' = (\mu_0', \mu_1', \ldots, \mu_n')$,  and $\nu' = (\nu_0', \nu_1', \ldots, \nu_n')$. 
Suppose that we can find $i$, $j$, and $k$ in $\{0,\ldots, n\}$ such that $i+j=k+n$ and $\lambda'_{i}+\mu'_{j}=\nu'_{k}$.  Let 
$\tilde{\lambda}'$, $\tilde{\mu}'$, and $\tilde{\nu}'$ be the characters obtained by removing the 
$i^{\text{th}}$, $j^{\text{th}}$, and $k^{\text{th}}$ coordinates
from $\lambda'$, $\mu'$, and $\nu'$ respectively, i.e., 
\begin{align*}
\tilde{\lambda}' := (\lambda_{0}', \ldots, \lambda_{i-1}', \lambda_{i+1}',& \ldots, \lambda_{n}'), \,\,
\tilde{\mu}' := (\mu_{0}', \ldots, \mu_{j-1}', \mu_{j+1}',\ldots, \mu_{n}'),  \\
\mbox{and}\,\,\,\tilde{\nu}' & := (\nu_{0}', \ldots, \nu_{k-1}', \nu_{k+1}',\ldots, \nu_{n}').
\end{align*}
One checks immediately that this is the Littlewood-Richardson version of Reduction Formula I from \cite{GH},
and hence that
$c_{\lambda', \mu'}^{\nu'} = c_{\tilde{\lambda}', \tilde{\mu}'}^{\tilde{\nu}'}$.

To find such $i$, $j$, and $k$ we use Theorem \ref{theorem1}. 
Let $w_1$, $w_2$, and $w_3$ be the corresponding
Weyl group elements, i.e. $\Phi_{w_3} = \Phi_{w_1} \sqcup \Phi_{w_2}$ and $w_3^{-1} \cdot \nu' = w_1^{-1} \cdot \lambda' + w_2^{-1} \cdot \mu'$.
Since $\cW$ is the symmetric group on $n+1$ elements, we may think of elements of $\cW$ as bijective functions on $\{0, 1, \ldots, n\}$. 
For any $w \in \cW$ define the displacement function $\delta_w$ by $\delta_w(i) = w(i) - i$. It is clear that $\delta_w = \rho - w^{-1} \rho = -w^{-1} \cdot 0$
and hence $w^{-1} \cdot \tau = w^{-1} \tau - \delta_w$. Now $\Phi_{w_3} = \Phi_{w_1} \sqcup \Phi_{w_2}$ implies 
$w_3^{-1} \cdot 0 = w_1^{-1} \cdot 0 + w_2^{-1} \cdot 0$ which by the previous calculation may be written as 
\begin{equation} \label{equation41}
\delta_{w_3} = \delta_{w_1} + \delta_{w_2}.
\end{equation}
The identity $w_3^{-1} \cdot \nu' = w_1^{-1} \cdot \lambda' + w_2^{-1} \cdot \mu'$ together with (\ref{equation41}) gives
\begin{equation} \label{equation42}
(\nu_{w_3(0)}', \nu_{w_3(1)}', \ldots, \nu_{w_3(n)}') =   (\lambda_{w_1(0)}', \lambda_{w_1(1)}', \ldots, \lambda_{w_1(n)}') + (\mu_{w_2(0)}', \mu_{w_2(1)}', \ldots, \mu_{w_2(n)}').
\end{equation}
A comparison of the last coordinates of (\ref{equation41}) 
yields
$w_3(n) = \delta_{w_3}(n) + n = \delta_{w_1}(n) + \delta_{w_2}(n) + n = w_1(n) + w_2(n) - n$ while the last coordinates
of (\ref{equation42}) yield
$\nu'_{w_3(n)} = \lambda'_{w_1(n)} + \mu'_{w_2(n)}$. Therefore, setting $i := w_1(n)$,
$j:= w_2(n)$, and $k:= w_3(n)$, the previous identities become 
$$
i + j = k + n, \quad \mbox{and} \quad \lambda'_i + \mu'_j = \nu'_k.
$$

The final step is to verify that this particular reduction preserves the property of being a cohomological component,
i.e, that $V({\tilde{\nu}'})$ is a cohomological component of $V({\tilde{\lambda}'})\otimes V({\tilde{\mu}'})$.
Define permutations of $\{0, 1, \ldots, n-1\}$ by 
$$\tilde{w}_1(s)  = \begin{cases} w_1(s) & {\text { if }} w_1(s) < i \\ w_1(s) -1 & {\text { if }} w_1(s) >i \end{cases},$$ 
\begin{equation}\label{equation60}
\tilde{w}_2(s)  = \begin{cases} w_2(s) & {\text { if }} w_2(s) < j \\ w_2(s) -1 & {\text { if }} w_2(s) >j \end{cases},
\end{equation}
$$\tilde{w}_3(s)  = \begin{cases} w_3(s) & {\text { if }} w_3(s) < k \\ w_3(s) -1  & {\text { if }} w_3(s) >k \end{cases}.$$
It is not difficult to check that \eqref{equation60} ensures that  
$\tilde{w}_3^{-1}\cdot\tilde{\nu}'= \tilde{w}_1^{-1}\cdot\tilde{\lambda}'+\tilde{w}_2^{-1}\cdot\tilde{\mu}'$
and  
$\Phi_{\tilde{w}_3}=\Phi_{\tilde{w}_1}\sqcup\Phi_{\tilde{w}_2}$.
Therefore, by Theorem \ref{theorem1},
$V(\tilde{\nu}')$ is a cohomological component of $V(\tilde{\lambda}') \otimes V(\tilde{\mu}')$; i.e.
we have proved the following claim.

\begin{claim}\label{claim5}
Let $G = GL_{n+1}$ and let $V(\nu')$ be a cohomological component of $V(\lambda') \otimes V(\mu')$. Then the triple 
$(\lambda', \mu', \nu')$ admits a reduction pattern and therefore $c_{\lambda', \mu'}^{\nu'}=1$.
\end{claim}

Note that it is not true in general that an application of Reduction Formula I 
results in a cohomological component.
Let us call a sequence $\cR_n, \cR_{n-1}, \ldots, \cR_1$ of reductions which can
be applied consecutively to a triple $(\lambda', \mu', \nu')$ a weak reduction pattern (i.e., we do not
require that each step results in a cohomological component).
There exist weak reduction patterns which are
not reduction patterns, but
if at least one of the characters $\lambda', \mu'$, or $\nu'$ is strictly dominant, then any weak reduction pattern is a reduction pattern. 
We do not know 
whether the existence of a weak reduction pattern for $(\lambda', \mu', \nu')$ always implies the existence of 
a reduction pattern for $(\lambda', \mu', \nu')$. 
In other words, 
we do not know the answer to the following question.

\begin{question} \label{question5}
Is is true that if a triple $(\lambda', \mu', \nu')$ of dominant characters admits a weak reduction pattern then $V(\nu')$ is a cohomological component of 
$V(\lambda') \otimes V(\mu')$?
\end{question}
Notice that the Conjecture implies a positive answer to Question \ref{question5}. Indeed, if $(\lambda', \mu', \nu')$ admits a weak reduction pattern
then $V(\nu')$ is a generalized PRV component of 
$V(\lambda') \otimes V(\mu')$
of multiplicity one. We expect, however, that Question \ref{question5} may be answered by a direct combinatorial
argument.

\section{Cohomological components are extreme components of the tensor product.}
In this section we continue the discussion of Problem \ref{problem2}.
Theorem \ref{theorem2} gives a partial solution  -- which would be completed by the Conjecture -- but here
we present properties of cohomological components related to the combinatorics of the Littlewood--Richardson cone.

Denote by $LR$  the Littlewood--Richardson cone, i.e., 
the rational convex cone generated by  $\{(\lambda', \mu', \nu') \in (\cP^+)^3 \, | \,c_{\lambda', \mu'}^{\nu'} \neq 0\}$.
Given $\lambda', \mu' \in \cP^+$, let $LR(\lambda',\mu')$  be the slice of the
Littlewood-Richardson cone $LR$ obtained by fixing the first two coordinates to be $\lambda'$ and $\mu'$ respectively.  
The slice $LR(\lambda',\mu')$ is a convex polytope; this follows, for example, from the solution of Horn's conjecture, 
see \cite{Kly}, \cite{KT1}, and \cite{KM}.
The solution of Horn's conjecture also implies that a point $\nu'$ is in $LR(\lambda',\mu')$ if and only if there is a positive
integer $m$ such that $m\nu' \in \cP^+$ and such that $V_{m\nu'}$ is a component of $V_{m\lambda'}\otimes V_{m\mu'}$.

Part of the original 
theorem of Parthasarathy--Ranga Rao--Varadarajan implies that the highest weight of the PRV component is a vertex of  
$LR(\lambda', \mu')$.
Indeed, they proved that $V(\lambda') \otimes V(\mu')$ is generated by the vector $v_{\lambda'} \otimes v_{w_0 \mu'}$ where
$v_{\lambda'}$ is the highest weight vector of $V(\lambda')$ and $v_{w_0 \mu'}$ is the lowest weight vector of $V(\mu')$. In particular, the 
character of $v_\lambda \otimes v_{w_0 \mu'}$ is contained in the support of any irreducible submodule of $V(\lambda') \otimes V(\mu')$.  Since this also
holds after scaling by an arbitrary positive integer $m$, this implies that this character is a vertex of  $LR(\lambda', \mu')$. We have the following generalization of this fact.

\begin{claim} \label{claim6}
If $V(\nu')$ is a cohomological component of $V(\lambda') \otimes V(\mu')$, then $\nu'$ is a vertex of $LR(\lambda', \mu')$.
\end{claim}
The following pictures illustrate this statement. We have drawn  the polytope $LR(\lambda', \mu')$
and the highest weights of the components of $V(\lambda') \otimes V(\mu')$
in the cases $G = SL_3$ and $V((3,5)) \otimes V((1,2))$,
$G=SO_5$ and $V(\rho) \otimes V(\rho)$, and $G = SL_3$ and $V((7,2)) \otimes V((1,3))$. The cohomological 
components are circled and the generalized PRV components which are not cohomological are marked with a square.

\bigskip

\ifthenelse{\boolean{showpicture}}{

\noindent
\begin{tabular}{ccc}
\\
\begin{pspicture}(-0.8,1.5)(3.8,6.5)
\rput(1.5,1.1){\scriptsize {$V((3,5))\otimes V((1,2))$}}
\rput(1.5,6.9){\scriptsize {$SL_3$}}
\SpecialCoor
\pscustom[fillstyle=solid,fillcolor=superlightgray,linecolor=gray]{
\wtmvto{1}{4}
\wtlnto{5}{2}
\wtlnto{6}{3}
\wtlnto{4}{7}
\wtlnto{2}{8}
\wtlnto{0}{6}
\wtlnto{1}{4}
}
\Agrid{-0.8}{1.5}{3.8}{6.5}
\wtcirc{0}{6}
\wtcirc{1}{4}
\wtcirc{1}{7}
\wtcirc{2}{5}
\wtcirc{2}{8}
\wtcirc{3}{3}
\wtcirc{3}{6}
\wtcirc{4}{4}
\wtcirc{4}{7}
\wtcirc{5}{2}
\wtcirc{5}{5}
\wtcirc{6}{3}
\bigwtcirc{4}{7}
\bigwtcirc{6}{3}
\bigwtcirc{2}{8}
\bigwtcirc{5}{2}
\bigwtcirc{0}{6}
\bigwtcirc{1}{4}
\end{pspicture}
&
\begin{pspicture}(-1,-0.5)(5,4.5)
\rput(2,4.9){\scriptsize {$SO_5$}}
\rput(2,-0.9){\scriptsize {$V(\rho)\otimes V(\rho)$}}
\SpecialCoor
\pscustom[fillstyle=solid,fillcolor=superlightgray,linecolor=gray]{
\bwtmvto{0}{0}
\bwtlnto{3}{0}
\bwtlnto{2}{2}
\bwtlnto{0}{4}
\bwtlnto{0}{0}
}
\Bgrid{-0.5}{-0.5}{4.5}{4.5}
\bwtcirc{0}{0}
\bwtcirc{1}{0}
\bwtcirc{2}{0}
\bwtcirc{3}{0}
\bwtcirc{2}{2}
\bwtcirc{0}{4}
\bwtcirc{0}{2}
\bwtcirc{1}{2}
\bigbwtcirc{0}{0}
\bigbwtcirc{3}{0}
\bigbwtcirc{2}{2}
\bigbwtcirc{0}{4}
\bigbwtsq{1}{0}
\bigbwtsq{0}{2}
\bigbwtsq{1}{2}
\end{pspicture}
&
\begin{pspicture}(0.7,1)(5.8,6)
\rput(3.25,6.4){\scriptsize {$SL_3$}}
\rput(3.25,0.6){\scriptsize {$V((7,2))\otimes V((1,3))$}}
\SpecialCoor
\pscustom[fillstyle=solid,fillcolor=superlightgray,linecolor=gray]{
\wtmvto{4}{1}
\wtlnto{6}{0}
\wtlnto{9}{0}
\wtlnto{10}{1}
\wtlnto{8}{5}
\wtlnto{6}{6}
\wtlnto{3}{3}
\wtlnto{4}{1}
}
\Agrid{0.7}{1}{5.8}{6}
\wtcirc{3}{3}
\wtcirc{4}{1}
\wtcirc{4}{4}
\wtcirc{5}{2}
\wtcirc{5}{5}
\wtcirc{6}{0}
\wtcirc{6}{3}
\wtcirc{6}{6}
\wtcirc{7}{1}
\wtcirc{7}{4}
\wtcirc{8}{2}
\wtcirc{8}{5}
\wtcirc{9}{0}
\wtcirc{9}{3}
\wtcirc{10}{1}
\bigwtcirc{8}{5}
\bigwtcirc{10}{1}
\bigwtcirc{6}{6}
\bigwtcirc{3}{3}
\bigwtcirc{4}{1}
\bigwtsq{8}{2}
\end{pspicture}
\\
\\
\\
\multicolumn{3}{c}{\sc Figure 1} \\
\end{tabular}
}{}


There are several observations we can make about the diagrams above. 
In the first case all generalized PRV components are cohomological, while in the other two only some of them are.
Furthermore, in the first two cases all vertices of $LR(\lambda', \mu')$ are cohomological, 
while in the last one two vertices of 
$LR(\lambda', \mu')$ are not 
(although they are not generalized PRV components either). 
Finally, the number of cohomological components varies from $4$ to $6$. 
These observations lead to the following natural questions.

\begin{question} \label{question7}
Is it true that a generalized PRV component of $V(\lambda') \otimes V(\mu')$ is cohomological if and only if it is a vertex of 
$LR(\lambda', \mu')$?
\end{question}

\begin{question} \label{question8}
How many cohomological components does $V(\lambda')  \otimes V(\mu')$ have? 
How many vertices does $LR(\lambda', \mu')$ have? 
\end{question}
As a first step towards answering Question \ref{question8} one may try to determine the possible number of
cohomological components or vertices. In studying these questions it makes sense to count
the generalized PRV components with multiplicities, i.e. we count $V(\sigma'(\lambda' + w' \mu'))$ and 
$V(\sigma''(\lambda' + w'' \mu'))$ with $w' \neq w''$ as different generalized PRV components even if
$\sigma'(\lambda' + w' \mu') = \sigma''(\lambda' + w'' \mu'')$. Thus the number of generalized PRV components of
$V(\lambda') \otimes V(\mu')$ is equal to the order of $\cW$. With this convention in mind we see that 
$V(\lambda') \otimes V(\mu')$ has at most $|\cW|$  and at least $2^n$ cohomological components, where $n$ is the
semisimple rank of $G$. Indeed, for any parabolic subalgebra $P$ of $G$ containing $B$, let $w_0(P)$ be the longest
element of $\cW$ which stabilizes the roots of $P$. The generalized PRV component corresponding to $w_0(P)$ is 
parabolically induced from the PRV component of the reductive part of $P$. A moment's thought shows that 
parabolic induction preserves the property of a component being cohomological. Since there are exactly $2^n$
parabolic subgroups of $G$ containing $B$ (including $B$ and $G$ themselves), we conclude that there are 
at least that many cohomological components. Both of these limits are achievable. 
Claim \ref{claim2}(ii) provides us with examples in which there are  $|\cW|$ cohomological
components. The next example shows that the number $2^n$ can also be achieved.

\begin{example} \label{example9}
Consider the tensor product $V(\rho) \otimes V(\rho)$. The first observation is that if $V(\nu')$ is a cohomological component
then $\nu' = \rho + w_3 \rho \in \cP^+$. The second observation is that $\rho + w \rho \in \cP^+$ with $w \in W$ implies that
$w = w_0(P)$ for some parabolic subgroup $P$. Both of these observations are easy exercises that we leave to the reader 
(and/or their students). Thus we obtain $2^n$ distinct cohomological components among the (not distinct) generalized PRV 
components. We can also modify this example to obtain $2^n$ cohomological components in a case when all generalized
PRV components are distinct. Indeed, it is enough to take the tensor product $V((N+1) \rho) \otimes V(N \rho)$ for a
large enough $N$. The cohomological components are still parabolically induced from PRV components and the generalized
PRV components are all distinct.
\end{example}

The tensor product considered in the example above seems to be more approachable than the general case while still
retaining a very interesting structure. It may be a good test for some of the questions we stated throughout the paper
and in particular for the Conjecture. We do not know whether the Conjecture holds for $V(\rho) \otimes V(\rho)$.

We end the paper with two more questions concerning the combinatorics of cohomological components. Unlike Question \ref{question8}
where we are interested in the cohomological components of a particular tensor product, the following questions 
ask about cohomological components ``at large".

\begin{question} \label{question1}
Describe all triples $(w_1, w_2, w_3)$ which satisfy $\Phi_{w_3} = \Phi_{w_1} \sqcup \Phi_{w_2}$. 
\end{question}
This question is intentionally stated in a very general and vague way. There are different combinatorial interpretations of $\Phi_{w_3} = \Phi_{w_1} \sqcup \Phi_{w_2}$, and
we have already discussed some of them in Section 2.  To explain what kind of an answer we are looking for we turn to Demazure's proof of the
Borel--Weil--Bott theorem, \cite{De}. The proof boils down to showing that the Weyl group ``acts" on the cohomology of line bundles on $X$ in a 
way compatible with the affine action on characters. The answer to Question \ref{question1} we are hoping for would construct a graph whose 
vertices are the triples under consideration and whose edges correspond to natural transformations from a triple to another triple.
This very general question is related to some of the questions we stated above.

A very concrete version of Question \ref{question1}
is to ask for the number of such triples. We do not know the answer, but here is a related interesting fact. 
If $G = SL_{n+1}$, the number of $n$--tuples $(w_1, w_2, \ldots, w_n)$ of
nontrivial elements of $\cW$ such that $\Delta^+ = \Phi_{w_1} \sqcup \Phi_{w_2} \sqcup \ldots \sqcup \Phi_{w_n}$ is exactly the $n^{\rm{th}}$
Catalan number $\frac{1}{n+1}\binom{2n}{n}$.

For the last question we fix a triple $(w_1, w_2, w_3)$ as in Question \ref{question1} and consider the rational cone
$C_{w_1, w_2}^{w_3}$ generated by the set 
$$
\{ (\lambda', \mu', \nu') \in (\cP^+)^3 \, | \, w_3^{-1} \nu' = w_1^{-1} \lambda' + w_2^{-1} \mu'\} \rule{5cm}{0cm}
$$
$$
\rule{4cm}{0cm}
=
\{ (\lambda', \mu', \nu') \in (\cP^+)^3 \, | \, w_3^{-1} \cdot \nu' = w_1^{-1} \cdot \lambda' + w_2^{-1} \cdot \mu'\},
$$
where the equality of sets above is due to the equation 
$w_3^{-1}\cdot 0 = w_1^{-1}\cdot 0 + w_2^{-1}\cdot 0$, implied by the condition $\Phi_{w_3}=\Phi_{w_1}\sqcup \Phi_{w_2}$.
The cone $C_{w_1, w_2}^{w_3}$ is not empty as it always contains the triple $(0,0,0)$.

\begin{question} \label{question2}  
What is the dimension of $C_{w_1, w_2}^{w_3}$? 
\end{question}

The cones $C_{w_1, w_2}^{w_3}$ are connected to many of the previous questions.
As an example,
if 
$C_{w_1, w_2}^{w_3}$ contains a point $(\lambda', \mu', \nu')$ with at least one strictly dominant entry
then $d_{w_1, w_2}^{w_3}=1$, giving an answer to Question \ref{question0} for the triple $(w_1, w_2, w_3)$.

Finally we  note that combinations of the open questions stated in the paper may also be posed. For example, we can combine 
Question \ref{question8} and Question \ref{question2} by asking whether the triple $(w_1, w_2, w_3)$ determines the
number of cohomological components and if so, how. 

\subsection*{Note added in proof.} After this paper was accepted for publication we obtained an almost complete
proof of the Conjecture.  More precisely, we have proved (see \cite{dr}) that the Conjecture holds





\begin{enumerate}
\item[1.] when $G$ is a simple classical group;
\item[2.] when $G$ is any semisimple group and at least one of the characters $\lambda'$, $\mu'$, or $\nu'$ is strictly dominant.
\end{enumerate}

\vskip1cm
\begin{tabular}{ll}\noindent
Address: &Department of Mathematics and Statistics, Queen's University\\
& Kingston, K7L 3N6, Canada \\
E-mail: & I. D.: {\tt dimitrov@mast.queensu.ca} \\
& M. R.:  {\tt mikeroth@mast.queensu.ca}
\end{tabular}

\end{document}